\documentclass[letterpaper, 10 pt, conference]{ieeeconf}%
\usepackage{graphics}
\usepackage{epsfig}
\usepackage{cite}
\usepackage{breakurl}
\usepackage{amsmath}
\usepackage{amssymb}
\usepackage{amsfonts}
\usepackage{graphicx}
\usepackage[normalem]{ulem}
\usepackage{color}
\usepackage{dsfont}
\usepackage[hyphens]{url}
\usepackage[hidelinks]{hyperref}%
\providecommand{\U}[1]{\protect\rule{.1in}{.1in}}
\IEEEoverridecommandlockouts
\newtheorem{innercustomdef}{Definition}
\newenvironment{customdef}[1]
  {\renewcommand\theinnercustomdef{#1}\innercustomdef}
  {\endinnercustomdef}
\newtheorem{innercustomass}{Assumption}
\newenvironment{customass}[1]
  {\renewcommand\theinnercustomass{#1}\innercustomass}
  {\endinnercustomass}

\allowdisplaybreaks[0]
\hypersetup{breaklinks=true}
\begin{document}

\title{{\LARGE \textbf{The Penetration Effect of Connected Automated Vehicles in Urban Traffic: An Energy Impact Study}}}
\author{Yue Zhang, and Christos G. Cassandras \thanks{ The work of Cassandras and Zhang is supported in part by NSF under grants ECCS-1509084, CNS-1645681, and IIP-1430145, 
by AFOSR under grant FA9550-15-1-0471, by DOE under grant DOE-46100, by MathWorks and by Bosch.} \thanks{Y. Zhang and C.G. Cassandras are with the Division of
Systems Engineering and Center for Information and Systems Engineering, Boston
University, Boston, MA 02215 USA (e-mail: joycez@bu.edu; cgc@bu.edu).} }
\maketitle

\begin{abstract}
Earlier work has established a decentralized framework of optimally
controlling connected and automated vehicles (CAVs) crossing an urban
intersection without using explicit traffic signaling. The proposed solution
is capable of minimizing energy consumption subject to a throughput
maximization requirement. In this paper, we
address the problem of optimally controlling CAVs under mixed traffic
conditions where both CAVs and human-driven vehicles (non-CAVs) travel on
the roads, so as to minimize energy consumption while guaranteeing safety constraints. The impact of
CAVs on overall energy consumption is also investigated under different traffic
scenarios. The benefit from CAV penetration (i.e., the fraction of CAVs relative to all vehicles) is validated through
simulation in MATLAB and VISSIM. The results indicate that the energy
efficiency improvement becomes more significant as the CAV penetration rate
increases, while the significance diminishes as the traffic  becomes heavier.

\end{abstract}

\thispagestyle{empty} \pagestyle{empty}


\section{Introduction}

To date, traffic light signaling is the prevailing method used for controlling
the traffic flow at urban intersections. Exploiting data-driven control and optimization approaches, recent research (e.g., see \cite{Fleck2015}) has enabled the adaptive adjustment of traffic light cycles, leading to reduced congestion. However, aside from the
obvious infrastructure cost, 
urban traffic lights can lead to more rear-end collisions and
reduced safety. These issues have motivated research efforts on new
approaches for signal-free intersection traffic control.

Connected and Automated Vehicles (CAVs) possess the potential to change the
transportation landscape by enabling users to better monitor
transportation network conditions and to improve traffic flow in terms of
reducing energy consumption, travel delays, accidents and greenhouse gas
emissions. One of the very early efforts was proposed in \cite{Levine1966},
where an optimal linear feedback regulator is designed to control a string of
vehicles. Dresner and Stone \cite{Dresner2004} proposed a reservation-based
scheme for automated vehicle intersection management, whereby a centralized
controller coordinates the vehicle crossing sequence based on the request
received from the vehicles. Numerous efforts based on reservation schemes
have been reported in the literature \cite{Dresner2008, DeLaFortelle2010,
Huang2012}. Reducing the travel delay and increasing the throughput of an
intersection is one desired goal to be achieved. Relevant efforts include \cite{Li2006,Yan2009,Zohdy2012, Zhu2015} which aim at minimizing the vehicle
travel time under collision-avoidance constraints. Lee and Park \cite{Lee2012}
focused on minimizing the overlap between vehicle positions. Miculescu and
Karaman \cite{Miculescu2014} used a polling-system to model the
intersection. 
A detailed discussion of the research in
the literature can be found in \cite{Rios-Torres}.

Our earlier work \cite{ZhangMalikopoulosCassandras2016} has established a
decentralized optimal control framework for coordinating online a continuous
flow of CAVs crossing an urban intersection without using explicit traffic
signaling. For each CAV, an energy minimization optimal control problem is
formulated where the time to cross the intersection is first determined through a
throughput maximization problem. We also established conditions under
which feasible solutions to the optimal control problem exist.

The benefits of CAV coordination and control on energy consumption have been established and quantified in recent literature  \cite{gilbert1976vehicle, hooker1988optimal,
hellstrom2010design, li2012minimum}. However, 
the integration of CAVs with conventional vehicles faces several challenges before their penetration rate (i.e., the fraction of CAVs relative to all vehicles in a transportation system) becomes significant. Thus, a critical question is that of determining the penetration effect of CAVs under mixed traffic conditions.  
Under such conditions, it is necessary to design control algorithms for CAVs and coordination policies that can accommodate both CAVs and conventional human-driven vehicles. Dresner and Stone \cite{dresner2007sharing} proposed a light model that can control the physical traffic lights as well as implementing a reservation-based control algorithm for autonomous vehicles while ensuring safety. Other efforts include using information from CAVs to better adapt the traffic light in a mixed traffic scenario (e.g., see \cite{guler2014using}).
In this paper, we address the problem of optimally controlling the CAVs
crossing an urban intersection in a mixed traffic scenario where both
CAVs and non-CAVs (conventional human-driven vehicles) travel on the roads. A decentralized optimal control
framework is presented whose solution yields the optimal
acceleration/deceleration so as to minimize the energy consumption subject to
a throughput maximizing requirement, while taking the interaction between
CAVs and non-CAVs into consideration.

The structure of the paper is as follows. In Section  \ref{sec:model}, we
review the model in \cite{ZhangMalikopoulosCassandras2016} and its
generalization in \cite{Malikopoulos2016}. In Section  \ref{sec:cav}, we
formulate the optimal control
problem for each CAV under mixed traffic conditions and present the
analytical solutions. In Section \ref{sec:mixed}, we introduce the non-CAV model and the approaches for collision avoidance among CAVs and non-CAVs to complete the establishment of the mixed traffic scenario. In Section
\ref{sec:sim}, we investigate the impact of CAV penetration on energy economy
under different traffic scenarios. We offer concluding remarks in Section
\ref{sec:con}.


\section{The Model}

\label{sec:model}

\begin{figure}[ptb]
\centering
\includegraphics[width=1\columnwidth]{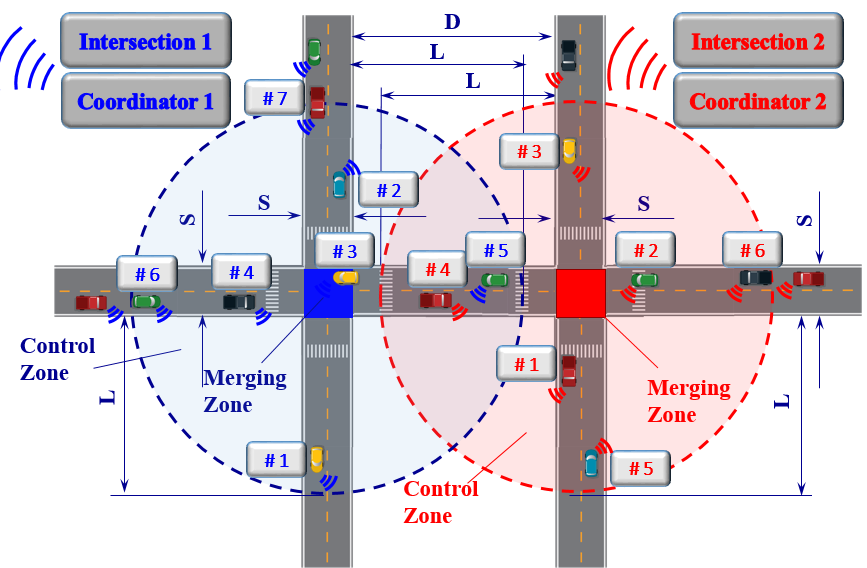} \caption{Connected Automated Vehicles (CAVs) crossing two adjacent intersections.}%
\label{fig:intersection}%
\end{figure}

We briefly review the model introduced in
\cite{ZhangMalikopoulosCassandras2016} where there are two intersections, 1
and 2, located within a distance $D$ (Fig. \ref{fig:intersection}). The region
at the center of each intersection, called \emph{Merging Zone} (MZ), is the
area of potential lateral CAV collision. Although it is not restrictive, this
is taken to be a square of side $S$. Each intersection has a \emph{Control
Zone} (CZ) and a coordinator that can communicate with the CAVs traveling
within it. The distance between the entry of the CZ and the entry of the MZ is
$L>S$, and it is assumed to be the same for all entry points to a given CZ.

Let $N_{z}(t)\in\mathbb{N}$ be the cumulative number of CAVs which have
entered the CZ of intersection $z$ at $t$ and formed a queue $\mathcal{N}_{z}(t)=\{1,\ldots,N_{z}(t)\}$
which designates the order in which these vehicles will be entering the MZ.
The way the queue is formed is not restrictive. When a CAV reaches the CZ of
intersection $z$, the coordinator assigns it an integer value $i=N_{z}(t)+1$.
If two or more CAVs enter a CZ at the same time, then the corresponding
coordinator selects randomly the first one to be assigned the value
$N_{z}(t)+1$. In the region between the exit point of a MZ and the entry point
of the subsequent CZ, the CAVs are assumed to cruise with the speed they had
when they exited that MZ.

For simplicity, we assume that each CAV is governed by second order dynamics:%
\begin{equation}
\dot{p}_{i}(t)=v_{i}(t)\text{, }~p_{i}(t_{i}^{0})=0\text{; }~\dot{v}%
_{i}(t)=u_{i}(t)\text{, }v_{i}(t_{i}^{0})\text{ given} \label{eq:model2}%
\end{equation}
where $p_{i}(t)$, $v_{i}(t)$, and $u_{i}(t)$
denote the position (i.e., travel distance since the entry of the CZ), speed
and acceleration/deceleration (i.e., control input) of each CAV $i$.
These dynamics are in force over an interval $[t_{i}^{0},t_{i}^{f}]$, where
$t_{i}^{0}$ and $t_{i}^{f}$ are the times that CAV $i$ enters the CZ
and exits the MZ of intersection $z$ respectively.

To ensure that the control input and vehicle speed are within a given
admissible range, the following constraints are imposed:
\begin{equation}%
\begin{split}
u_{i,min}  &  \leq u_{i}(t) \leq u_{i,max},\quad\text{and}\\
0  &  \leq v_{min} \leq v_{i}(t) \leq v_{max},\quad\forall t\in\lbrack
t_{i}^{0},t_{i}^{m}],
\end{split}
\label{speed_accel constraints}%
\end{equation}
where $t_{i}^{m}$ is the time that CAV $i$ enters the MZ.

\begin{customdef}{1}
Depending on its physical location inside the CZ, CAV $j\in\mathcal{N}%
_{z}(t), j\neq i$ belongs to only one of the following four subsets of $\mathcal{N}%
_{z}(t)$ with respect to CAV $i$: 1) $\mathcal{R}^{z}_{i}(t)$ contains all
CAVs traveling on the same road as $i$ and towards the same direction but on
different lanes, 2) $\mathcal{L}^{z}_{i}(t)$ contains all CAVs traveling on
the same road and lane as vehicle $i$ (e.g., $\mathcal{L}^{1}_{6}(t)$ contains
CAV \#4 in Fig. \ref{fig:intersection}), 3) $\mathcal{C}^{z}_{i}(t)$ contains
all CAVs traveling on different roads from $i$ and having destinations that
can cause collision at the MZ (e.g., $\mathcal{C}^{1}_{7}(t)$ contains CAV \#6
in Fig. \ref{fig:intersection}), and 4) $\mathcal{O}^{z}_{i}(t)$ contains all
CAVs traveling on the same road as $i$ and opposite destinations that cannot,
however, cause collision at the MZ (e.g., $\mathcal{O}^{1}_{4}(t)$ contains
CAV \#3 in Fig. \ref{fig:intersection}). \label{def:2}
\end{customdef}

To ensure the absence of any rear-end collision throughout the CZ, we impose
the \emph{rear-end safety constraint}:
\begin{equation}
s_{i}(t)=p_{k}(t)-p_{i}(t)\geq\delta,~\forall t\in\lbrack t_{i}^{0},t_{i}%
^{m}],~k\in\mathcal{L}_{i}^{z}(t) \label{rearend}%
\end{equation}
where 
$k$ is the CAV physically ahead of $i$ on the same lane, $s_i(t)$ is the inter-vehicle distance between $i$ and $k$, and $\delta$ is the \emph{minimal safety following distance} allowable.

A lateral collision involving CAV $i$ may occur only if some CAV $j\neq i$
belongs to $\mathcal{C}^{z}_{i}(t)$. This leads to the following definition:

\begin{customdef}{2}
For each CAV $i\in\mathcal {N}_{z}(t)$, we define the set $\Gamma_{i}$ that
includes all time instants when a lateral collision involving CAV $i$ is
possible: $\Gamma_{i}\triangleq\Big\{t~|~t\in\lbrack t_{i}^{m},t_{i}^{f}]\Big\}$.
\end{customdef}

Consequently, to avoid a lateral collision for any two vehicles $i,j\in
\mathcal{N}_{z}(t)$ on different roads, the following constraint should hold
\begin{equation}
\Gamma_{i}\cap\Gamma_{j}=\varnothing,\text{ \  }\forall t\in\lbrack
t_{i}^{m},t_{i}^{f}]\text{, \ }j\in\mathcal{C}^{z}_{i}(t). \label{eq:lateral}%
\end{equation}

As part of safety considerations, we impose the following assumption (which
may be relaxed if necessary):


\begin{customass}{1}
For CAV $i$, none of the constraints in (\ref{speed_accel constraints}%
)-(\ref{rearend}) is active at $t_{i}^{0}$. \label{ass:feas}
\end{customass}

\begin{customass}{2}
The speed of the CAVs inside the MZ is constant, i.e., $v_{i}(t)=v_{i}%
(t_{i}^{m})=v_{i}(t_{i}^{f})$, $\forall t\in\lbrack t_{i}^{m},t_{i}^{f}]$.
This implies that $t_{i}^{f}=t_{i}^{m}+\frac{S}{v_{i}(t_{i}^{m})}$. \label{ass:constant} 
\end{customass}


\begin{customass}{3}
Each CAV $i$ has proximity sensors and can measure local information without
errors or delays. \label{ass:sensor}
\end{customass}

The objective of each CAV is to derive an optimal acceleration/deceleration profile,
in terms of minimizing energy consumption, inside the CZ (i.e., over the time
interval $[t_{i}^{0},t_{i}^{m}]$) while avoiding congestion between the two
intersections. Since the coordinator is not involved in any decision making process on the
vehicle control, we can formulate $N_{1}(t)$ and $N_{2}(t)$ decentralized
tractable problems for intersection 1 and 2 respectively that can be solved
online. The terminal times for CAVs entering the MZ can be obtained as the
solutions to a throughput maximization problem formulated in
\cite{Malikopoulos2016} subject to rear-end and lateral collision avoidance
constraints inside the MZ. As shown in \cite{Malikopoulos2016}, the terminal time of CAV $i$ (i.e.,
$t_{i}^{m}$) can be recursively determined through%
\begin{equation}
t_{i}^{m^{\ast}}=\left\{
\begin{array}
[c]{ll}%
t_{1}^{m^{\ast}} & \mbox{if $i=1$}\\
\text{max }\{t_{i-1}^{m^{\ast}},t_{k}^{m^{\ast}}+\frac{\delta}{v_{k}^{m}%
},t_{i}^{c}\} & \text{if }i-1\in\mathcal{R}^{z}_{i}\cup\mathcal{O}^{z}_{i}\\
\text{max }\{t_{i-1}^{m^{\ast}}+\frac{\delta}{v_{i-1}^{m}},t_{i}^{c}\} &
\mbox{if $i-1\in\mathcal{L}^z_{i}$}\\
\text{max }\{t_{i-1}^{m^{\ast}}+\frac{S}{v_{i-1}^{m}},t_{i}^{c}\} &
\mbox{if $i-1\in\mathcal{C}^z_{i}$}
\end{array}
\right.  \label{def:tf}%
\end{equation}
where $t_{i}^{c}=t_{i}^{1}\mathds{1}_{v_{i}^{m}=v_{max}}+t_{i}^{2}%
(1-\mathds{1}_{v_{i}^{m}=v_{max}})$ and
\begin{align}
t_{i}^{1}  &  =t_{i}^{0}+\frac{L}{v_{max}}+\frac{(v_{max}-v_{i}^{0})^{2}%
}{2u_{i,max}v_{max}},\nonumber\\
t_{i}^{2}  &  =t_{i}^{0}+\frac{[2Lu_{i,max}+(v_{i}^{0})^{2}]^{1/2}-v_{i}^{0}%
}{u_{i,max}}.\nonumber
\end{align}
Here, $t_i^c$ is a lower bound of $t_i^m$ regardless of the solution of the throughput maximization problem.

The conditions under which the rear-end safety constraint in \eqref{rearend} does
not become active inside the CZ are provided in \cite{Zhang2016}.

\section{Optimal Control of CAVs in Mixed Traffic}

\label{sec:cav}

\begin{figure}[ptb]
\centering
\includegraphics[width=0.9\columnwidth]{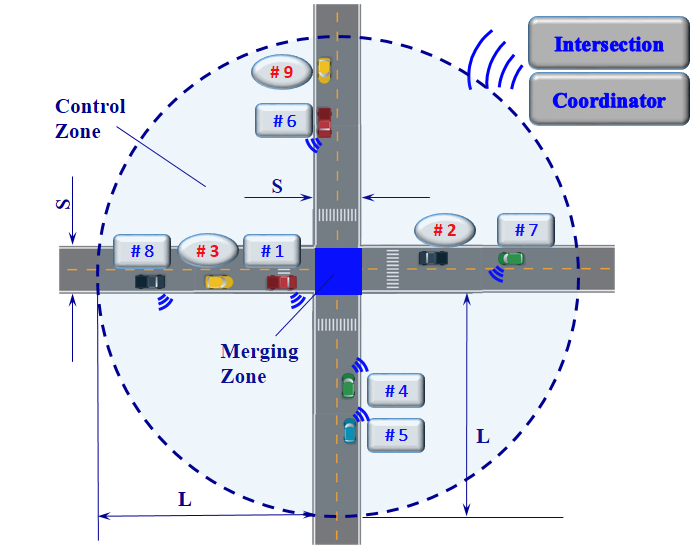} \caption{Connected
Automated Vehicles (blue labels) and non-CAVs (red labels) crossing an
urban intersection.}%
\label{fig:mixed}%
\end{figure}

We consider the mixed traffic scenario (Fig. \ref{fig:mixed}) where both CAVs
and non-CAVs (conventional human-driven vehicles) travel on the roads.
The first major issue to be addressed is modeling the interaction between CAVs
and non-CAVs where we assume that the latter do not possess the capability to communicate with other vehicles.

Regarding a CAV, there are two modes that it can be in: $(i)$ \textit{Free
Driving }(FD mode) when it is not constrained by a non-CAV that precedes it.
$(ii)$ \textit{Adaptive Following }(AF mode) when it follows a preceding
non-CAV while adaptively maintaining a safe following distance from it. CAVs
switch from the FD mode to the AF mode as soon as the inter-vehicle distance
falls below a certain threshold.

\subsection{Optimal Control for Free Driving (FD) Mode}

\label{oc_free} In this mode, the objective of each CAV is to derive an
optimal acceleration/deceleration profile, in terms of minimizing energy consumption, inside the
CZ, that is,
\begin{gather}
\min_{u_{i}\in U_{i}}\frac{1}{2}\int_{t_{i}^{0}}^{t_{i}^{m}}K_{i}\cdot
u_{i}^{2}(t)~dt\nonumber\\
\text{subject to}:\eqref{eq:model2},(\ref{speed_accel constraints}),t_{i}%
^{m},\text{ }p_{i}(t_{i}^{0})=0\text{, }p_{i}(t_{i}^{m}%
)=L,\label{eq:decentral}\\
\text{given }t_{i}^{0}\text{, }v_{i}(t_{i}^{0}),\nonumber
\end{gather}
where $K_{i}$ is a factor to capture CAV diversity (for simplicity, $K_{i}=1$
in the remainder of this paper). Note that this formulation does not include
the rear-end safety constraint in the CZ in \eqref{rearend}; we will return to this issue in what follows. On the other hand, the rear-end and lateral collision avoidance inside the MZ can be implicitly ensured by $t_i^m$.

An analytical solution of problem \eqref{eq:decentral} may be obtained through
a Hamiltonian analysis found in \cite{Malikopoulos2016}. Assuming that all
constraints are satisfied upon entering the CZ and that they remain inactive
throughout $[t_{i}^{0},t_{i}^{m}]$, the optimal control input
(acceleration/deceleration) over $t\in\lbrack t_{i}^{0},t_{i}^{m}]$ is given
by
\begin{equation}
u_{i}^{\ast}(t)=a_{i}t+b_{i} \label{eq:20}%
\end{equation}
where $a_{i}$ and $b_{i}$ are constants of integration. Using (\ref{eq:20}) in the CAV
dynamics \eqref{eq:model2}, the optimal speed and position are obtained:
\begin{equation}
v_{i}^{\ast}(t)=\frac{1}{2}a_{i}t^{2}+b_{i}t+c_{i} \label{eq:21}%
\end{equation}%
\begin{equation}
p_{i}^{\ast}(t)=\frac{1}{6}a_{i}t^{3}+\frac{1}{2}b_{i}t^{2}+c_{i}t+d_{i},
\label{eq:22}%
\end{equation}
where $c_{i}$ and $d_{i}$ are constants of integration. The coefficients
$a_{i}$, $b_{i}$, $c_{i}$, $d_{i}$ can be obtained given initial and terminal
conditions as follows:%

\begin{equation}
\left[
\begin{array}
[c]{cccc}%
\frac{1}{6}(t_{i}^{0})^{3} & \frac{1}{2}(t_{i}^{0})^{2} & (t_{i}^{0}) & 1\\
\frac{1}{2}(t_{i}^{0})^{2} & (t_{i}^{0}) & 1 & 0\\
\frac{1}{6}(t_{i}^{m})^{3} & \frac{1}{2}(t_{i}^{m})^{2} & t_{i}^{m} & 1\\
-t_{i}^{m} & -1 & 0 & 0
\end{array}
\right]  .\left[
\begin{array}
[c]{c}%
a_{i}\\
b_{i}\\
c_{i}\\
d_{i}%
\end{array}
\right]  =\left[
\begin{array}
[c]{c}%
p_{i}(t_{i}^{0})\\
v_{i}(t_{i}^{0})\\
p_{i}(t_{i}^{m})\\
\lambda_{i}^{v}(t_{i}^{m})
\end{array}
\right]  \label{eq:23}%
\end{equation}
Note that the analytical solution \eqref{eq:20} is valid while none of the
constraints becomes active for $t\in\lbrack t_{i}^{0},t_{i}^{m}]$. Otherwise,
the optimal solution should be modified considering the constraints, as
discussed in \cite{Malikopoulos2016}. Recall that the constraint
\eqref{rearend} is not included in \eqref{eq:decentral}. To address this
constraint, conditions under which the CAV is able to maintain feasibility in
terms of satisfying \eqref{rearend} over $t\in\lbrack t_{i}^{0},t_{i}^{m}]$
are derived in \cite{Zhang2016} along with an explicit mechanism to enforce
them prior to entering the CZ.

\subsection{Optimal Control for Adaptive Following (AF) Mode}

\label{oc_following}

\begin{figure}[ptb]
\centering
\includegraphics[width=1\columnwidth]{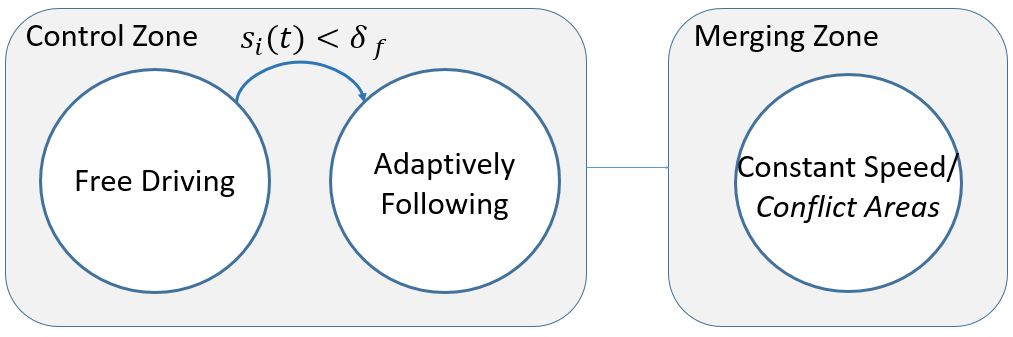} \caption{Modeling
approach for CAVs.}%
\label{CAVmodel}%
\vspace{-3mm}
\end{figure}

When the preceding vehicle for CAV $i$ is also a CAV, the feasibility of the
optimal solution (\ref{eq:20})-(\ref{eq:23}) can be enforced through an
appropriately designed \textit{Feasibility Enforcement Zone} that precedes the
CZ as described in \cite{Zhang2016}. Otherwise, when the inter-vehicle distance
$s_{i}(t)$ between CAV $i$ and the preceding vehicle $k$ falls below a certain
threshold $\delta_{f}$ at $t_{i}^{1}$, CAV $i$ transitions from the FD to the
AF mode (Fig. \ref{CAVmodel}). Since CAV $i$ cannot communicate with the
preceding non-CAV, it simply assumes a constant speed for the non-CAV. 

In this mode, the objective of each CAV is to derive an optimal
acceleration/deceleration profile so as to minimize energy consumption, while
maintaining the \textit{minimum safety following distance} $\delta$ with the
preceding non-CAV, that is,
\begin{gather}
\min_{u_{i}\in U_{i}}\frac{1}{2}\int_{t_{i}^{1}}^{t_{i}^{m}}[w_{u}\cdot
u_{i}^{2}(t)+w_{s}\cdot(s_{i}(t)-\delta)^{2}]~dt\nonumber\\
\text{subject to}:\eqref{eq:model2},(\ref{speed_accel constraints}),t_{i}%
^{m},\text{ }p_{i}(t_{i}^{m})=L,\label{eq:decentral2}\\
\text{and given }t_{i}^{1}\text{, }v_{i}(t_{i}^{1}),p_{i}(t_{i}^{1}),\nonumber
\end{gather}
where $w_{u}$ and $w_{s}$ are weights applied to the objective function, which
allow trading off energy consumption minimization against maintaining the
safety following distance.

The analytical solution of problem \eqref{eq:decentral2} may be obtained
through a Hamiltonian analysis similar to that in \cite{ZhangMalikopoulosCassandras2016} and \cite{Malikopoulos2016}. Assuming that all constraints are satisfied at
$t_{i}^{1}$ and that they remain inactive throughout $[t_{i}^{1},t_{i}^{m}]$,
the optimal control input (acceleration/deceleration) over $t\in\lbrack
t_{i}^{1},t_{i}^{m}]$ can be determined as follows. Defining $w=\frac
{\sqrt{4w_{u}\cdot w_{s}}}{2w_{u}}$ and $\alpha=\sqrt{\frac{w}{2}}$, the
optimal control can be obtained as
\begin{equation}
\begin{aligned} u_{i}^{*}(t) & = - 2a_{i}  \alpha^2 e^{\alpha t} sin(\alpha t) + 2b_{i}  \alpha^2 e^{-\alpha t} sin(\alpha t) \\ & + 2 c_{i}  \alpha^2 e^{\alpha t} cos(\alpha t) - 2d_{i}  \alpha^2 e^{-\alpha t} cos(\alpha t). \label{u_acc3} \end{aligned}
\end{equation}
The optimal speed and position can be obtained according to \eqref{eq:model2}:%
\begin{equation}
\begin{aligned} v_i^*(t) & = a_i \cdot \alpha e^{\alpha t}(cos(\alpha t) -sin(\alpha t)) \\ & - b_i \cdot \alpha e^{-\alpha t} (cos(\alpha t) + sin(\alpha t)) \\ & + c_i \cdot \alpha e^{\alpha t} (cos(\alpha t) + sin(\alpha t)) \\ & + d_i \cdot \alpha e^{-\alpha t} (cos(\alpha t) - sin(\alpha t)) \\ \end{aligned}
\end{equation}
\begin{equation}
\begin{aligned} p_i^*(t) & = a_{i} \cdot e^{\alpha t} cos(\alpha t) + b_{i}\cdot e^{-\alpha t} cos(\alpha t) \\ & + c_{i} \cdot e^{\alpha t} sin(\alpha t) + d_{i} \cdot e^{-\alpha t} sin(\alpha t) \\ & + p_{i-1}(t_{i}^{1}) + v_{i-1}(t_{i}^{1}) (t-t_{i}^{1}) - \delta, \end{aligned}
\end{equation}
where $a_{i}$, $b_{i}$, $c_{i}$ and $d_{i}$ are constants of integration,
which can be obtained in a similar way as \eqref{eq:23}.

\subsection{Terminal Conditions}

\label{oc_tim}

Under mixed traffic conditions, the recursive terminal time structure in
(\ref{def:tf}) derived in \cite{Malikopoulos2016} can no longer be applied
since non-CAVs are not controlled to follow the order imposed by the
queueing structure discussed in Sec. II. To determine the terminal
conditions for CAV $i$, there are two different cases to consider: $(i)$
vehicle $i-1$ is a CAV, and $(ii)$ vehicle $i-1$ is a non-CAV. If vehicle
$i-1$ is a CAV, then the terminal time for CAV $i$ can be recursively
determined through CAV $i-1$ as in \cite{Malikopoulos2016}; if vehicle $i-1$
happens to be a non-CAV, then the terminal time for CAV $i$ is determined by
estimating the terminal time of vehicle $i-1$. In particular, at time
$t_{i}^{0}$, vehicle $i-1$ is at position $p_{i-1}(t_{i}^{0})$ with speed
$v_{i-1}(t_{i}^{0})$, which can be measured by CAV $i$ through on-board
sensors or through the coordinator. As CAV $i$ cannot communicate with non-CAV
$i-1$, it simply assumes a constant speed for $i-1$, i.e., $v_{i-1}%
(t)=v_{i-1}(t_{i}^{0})$ for $t\in\lbrack t_{i}^{0},t_{i}^{m}]$. Denoting the estimated terminal time for vehicle $i-1$ as $\hat{t}_{i-1}^{m}$, we have%

\begin{equation}
\begin{aligned} \hat{t}_{i-1}^m &= t_i^0 + \frac{L-p_{i-1}(t_i^0)}{v_{i-1}(t_i^0)}, \end{aligned}
\end{equation}
based on which, CAV $i$ can determine its own terminal time for entering the
MZ using (\ref{def:tf}). Note that the estimation may need re-evaluation in the case that non-CAV $i-1$ changes speed.

\subsection{Simulation Example}

\label{oc_sim}

\begin{figure}[tbh]
\centering
\includegraphics[width=1\columnwidth]{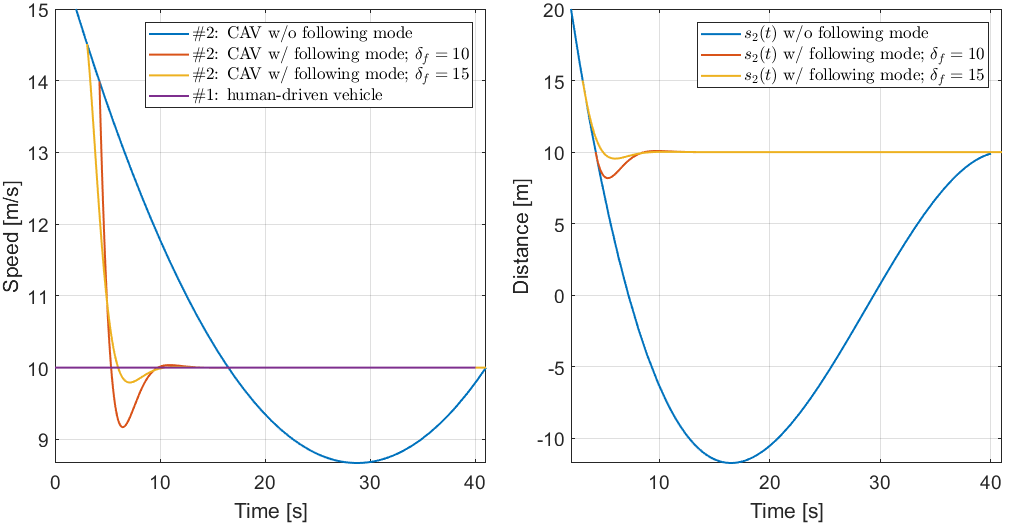}
\caption{Speed $v_{i}(t)$ and inter-vehicle distance $s_{2}(t)$ trajectories
w/ and w/o the optimal control for adaptively following.}%
\label{oc_exp}%
\end{figure}

The proposed optimal control framework for CAVs is illustrated through the
following simulation example, where the length of the CZ is $L=400$m. Vehicle
\#1 is assumed to be a non-CAV entering the CZ at $t_{1}^{0}=0$ and cruising
at its initial speed $v_{1}^{0}=10$m/s. CAV \#2 enters the same lane as
vehicle \#1 at $t_{2}^{0}=2$s with an initial speed $v_{2}^{0}=15$m/s. The
\textit{minimum safety following distance} is set to $\delta=10$m. With the
weights $w_{u}$ and $w_{s}$ both set to be 1, the optimal speed trajectory of
CAV \#2 (i.e., $v_2(t)$) and the inter-vehicle distance between vehicle \#1 and CAV \#2 (i.e., $s_2(t)$) with and without the
AF mode are shown in Fig. \ref{oc_exp}.


Without a transition to the AF mode, CAV \#2 violates the rear-end collision constraint (blue
curve in Fig. \ref{oc_exp}) and the distance $s_{2}(t)$ falls below
$\delta=10$m for $t>4.22$s. With the AF mode in force (red and yellow curves),
when the distance reaches the threshold, i.e., $s_{2}(t)=\delta_{f}$, CAV \#2
first decelerates so as to reach a much lower speed than vehicle \#1, and then
seeks to keep the distance as close to $\delta=10$m as possible. Note that the
process of adaptively following implicitly forces CAV \#2 to maintain the same
speed as vehicle \#1.

For the scenarios using the AF mode, the threshold for entering this mode is
set to either $\delta_{f}=10$m (red curve) or $\delta_{f}=15$m (yellow curve).
The energy consumption, given by a polynomial function of speed and
acceleration, is obtained as 0.0156 and 0.0143, respectively. The energy cost
reduction results from the fact that CAV \#2 enters the AF mode earlier given
$\delta_{f}=15$m, hence, it does not need to decelerate as hard as in the case
with $\delta_{f}=10$m. The approaching process becomes smoother, which leads
to lower energy consumption.



\section{Modeling Methodology for non-CAVs}

\label{sec:mixed}

There are two major issues that need to be addressed in terms of modeling
non-CAVs (i.e., human-driven vehicles): $(i)$ modeling the car-following
behavior, and $(ii)$ designing a collision avoidance approach inside the MZ
without explicit traffic signaling.

\subsection{The Wiedemann Approach}

\label{noncav:a} In this paper, we apply the Wiedemann model
\cite{wiedemann1974simulation}, the default approach adopted by VISSIM, the
transportation system simulator, to model car-following behavior. The basic
idea of the Wiedemann model is that a non-CAV can be in one of the
following  four driving modes:

\begin{itemize}
\item \textit{Free driving}: No observable influence by any preceding vehicle.
In this mode, the driver seeks to reach and (approximately) maintain a desired speed.

\item \textit{Approaching}: The driver adapts his/her own speed to the speed
of a preceding vehicle. This is done by decelerating until the speed
difference becomes zero when some desired safe distance from the vehicle being
followed is reached.

\item \textit{Following}: The driver follows the preceding vehicle while
trying to (approximately) maintain a safe distance from the vehicle being
followed. 

\item \textit{Braking}: The driver applies the brake to decelerate if the
distance from the preceding vehicle falls below the desired safety level. This
scenario may occur if the preceding vehicle changes speed abruptly.
\end{itemize}

For each mode, there are several associated parameters that model specific
car-following behavior. For instance, a \textquotedblleft\textit{Smooth
Closeup Behavior\textquotedblright} parameter can be enabled to model less
aggressive approaching behavior \cite{wiedemann1974simulation}. 

\subsection{Conflict Areas}

\label{noncav:b} As non-CAVs may not follow the prescribed order in the
queueing structure specified in Sec. II, lateral collisions may occur
inside the MZ when no traffic lights are present. There are several ways used
in VISSIM to model non-signalized intersections for non-CAVs, by defining
\textit{Priority Rules}, \textit{Conflict Areas}, and \textit{Stop Sign
Control}. Among these techniques, \textit{Conflict Areas} provide modeling
ease and more intelligent behavior and will be adopted in the sequel to ensure
the absence of lateral collisions in the MZ.

\begin{figure}[ptb]
\centering
\includegraphics[width=1\columnwidth]{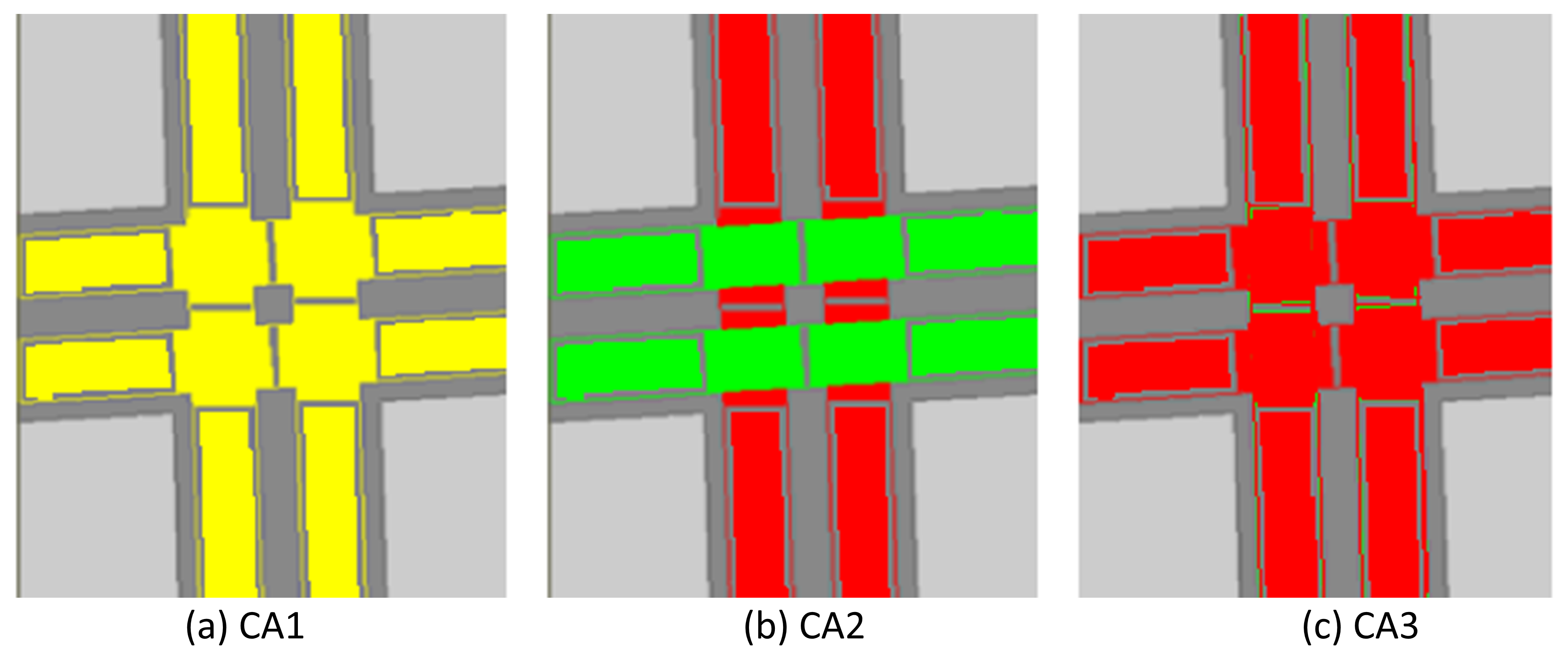}
\caption{Different states of conflict areas.}%
\label{conflictarea}%
\vspace{-4mm}
\end{figure}

If a CAV enters the MZ at the designated terminal time $t_{i}^{m}$ while a
non-CAV is present inside the MZ, the CAV simply forgoes the constant speed
assumption (\emph{Assumption 2}) and follows the \textit{Conflict Areas} rule
so as to avoid lateral collision. As shown in Fig. \ref{conflictarea}, there
are three options for defining conflict areas. Fig. \ref{conflictarea}(a)
indicates a passive conflict area, i.e., all vehicles are uncontrolled in
terms of lateral collision avoidance (CA1). Fig. \ref{conflictarea}(b) shows a
partially controlled conflict area, where the vehicles on the main road
(green) are uncontrolled because they have priority to cross the MZ, while
vehicles on the minor road (red) are controlled because they have to yield
(CA2). Fig. \ref{conflictarea}(c) shows a fully controlled conflict area. As
there is no right of way, vehicles on both roads are under control (CA3). In
this paper, we use the second conflict avoidance rule (i.e., CA2) for lateral collision
avoidance in mixed traffic. The third option is not appropriate
as it may lead to traffic deadlock.

\begin{figure}[tbh]
\centering
\includegraphics[width=1\columnwidth]{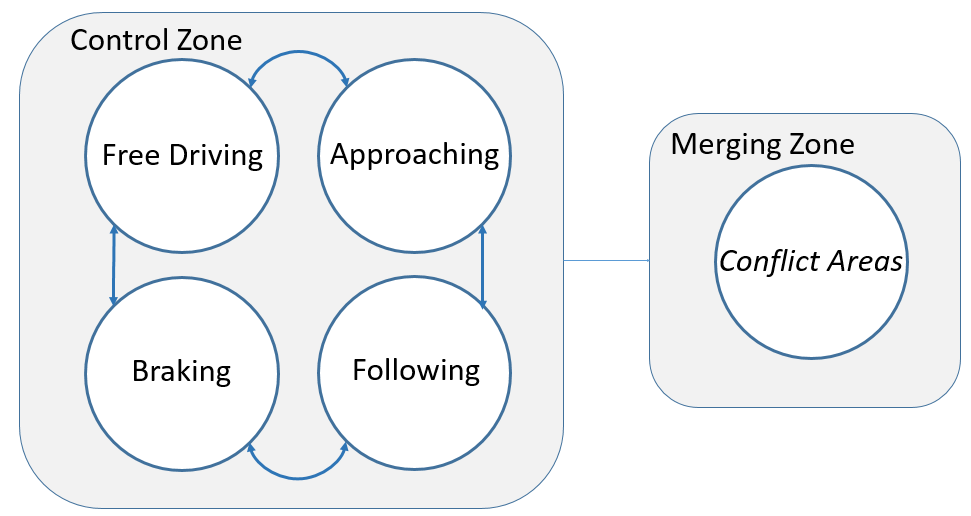} \caption{The
Wiedemann model for non-CAVs.}%
\label{nonCAVmodel}
\vspace{-3mm}
\end{figure}

The modeling approach for non-CAVs is summarized in Fig. \ref{nonCAVmodel},
where blue arrows indicate mode (state) transitions. The driver switches from
one mode to another as soon as a certain \emph{threshold} is reached, usually
expressed as a combination of speed difference and inter-vehicle distance.

\section{Energy Impact of CAV Penetration Under Different Traffic Scenarios}

\label{sec:sim}

The energy impact study is carried out through a combination of MATLAB and
VISSIM simulations. We consider a group of CAVs and non-CAVs crossing a single
urban intersection, where the length of the CZ is $L=400$m and the length of
the MZ is $S=30$m. For each direction, only one lane is considered. The
\textit{minimum safe following distance} is set to $\delta=10$m and the
threshold for entering the AF mode is $\delta_{f}=10$m. The weights $w_{u}$
and $w_{s}$ in \eqref{eq:decentral2} are both set to 1. The vehicle arrivals are assumed to be given by
a Poisson process and the initial speeds are uniformly distributed over
$[10.9,11.1]$m/s.

\subsection{Energy Impact of CAV Penetration}

\label{sim:a} 
\begin{figure}[tbh]
\centering
\includegraphics[width=1\columnwidth]{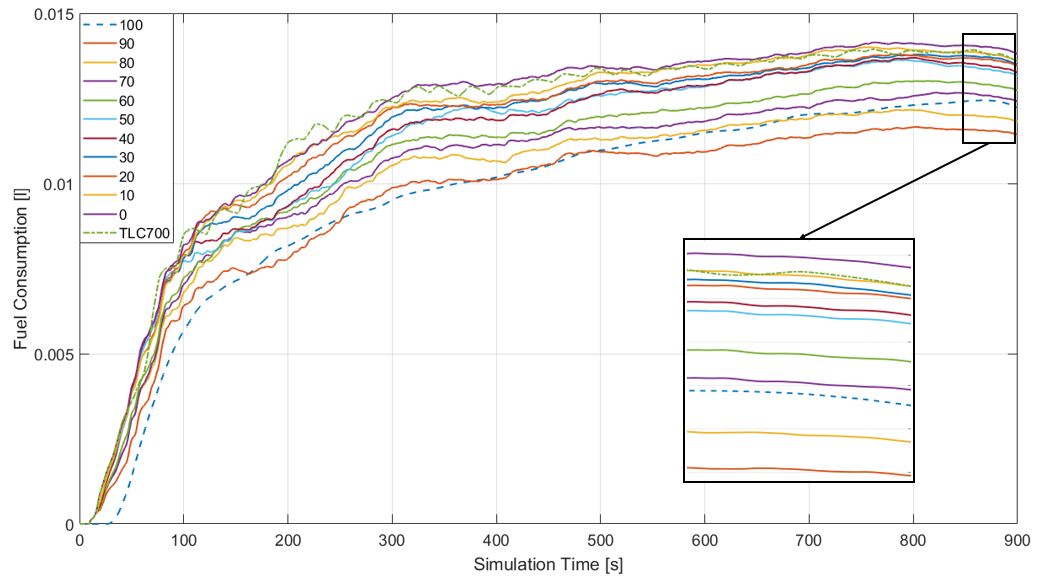} \caption{Energy
consumption per second with respect to different CAV penetration rates given
traffic flow rate set to $\lambda=$ 700 veh/(hour$\cdot$lane).}%
\label{700}%
\end{figure}

We first compare the energy impact over different CAV
penetration rates. Note that with 100\% CAV penetration, all CAVs proceed
according to the optimal trajectories determined in Sec. III and reach the
MZ at the designated terminal times. However, for cases with less than
100\% CAV penetration, we adopt the non-signalized collision avoidance rule (i.e., \textit{Conflict Areas}) in mixed traffic, as non-CAVs may not follow the prescribed order in
the queueing structure specified in Sec. II. The energy consumption
with respect to different CAV penetration rates given the traffic flow rate set
to $\lambda=700$ veh/(hour$\cdot$lane) is shown in Fig. \ref{700}. It can be seen that as the CAV penetration
rate increases, the energy consumption decreases, which validates the
efficiency of CAV penetration in terms of improving energy economy.

\begin{figure}[tbh]
\centering
\includegraphics[width=1\columnwidth]{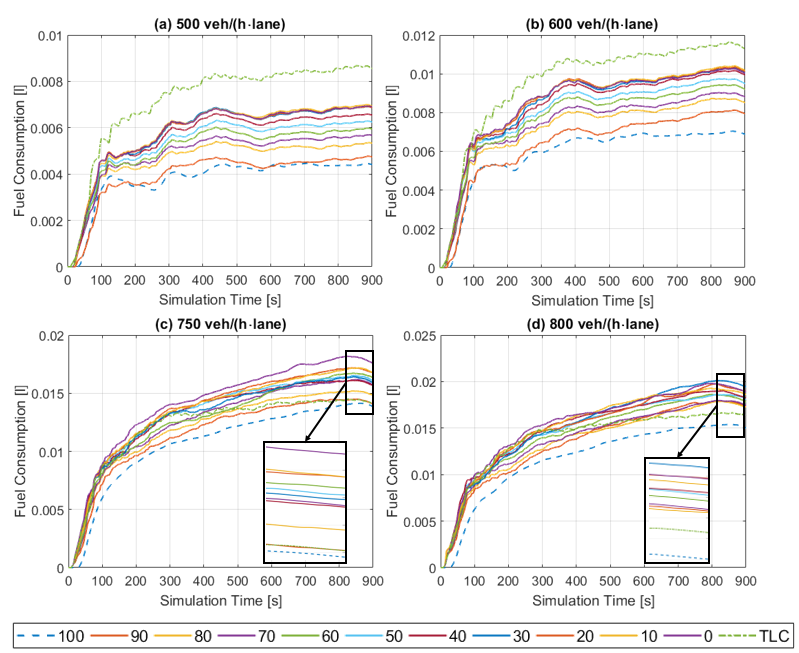} \caption{Energy
consumption per second with respect to different CAV penetration rates given
different traffic flow rates: (a) 500, (b) 600, (c) 750, (d) 800 veh/(hour$\cdot
$lane).}%
\label{678}%
\vspace{-3mm}
\end{figure}

Observe in Fig. \ref{700} that with no CAVs (0\% penetration rate), the
\textit{Conflict Areas} cannot outperform the traffic light control case
(indicated by TLC); however, with as little as 10\% CAV penetration rate, TLC is 
outperformed. This leads to the conclusion that approximately 10\% of
vehicles should be CAVs before energy consumption performance can exceed that
of TLC. However, this value clearly depends on traffic flow rates. 
Note that in Fig. \ref{678}(c) where the traffic flow rate is set to $\lambda=$ 750 veh/(hour$\cdot
$lane), we need 90\% of the vehicles to be CAVs in order to match the energy
performance under TLC; in Fig. \ref{678}(a) and \ref{678}(b) where the traffic flow rates are set to $\lambda = 500$ and $600$ veh/(hour$\cdot
$lane) respectively, the cases with no CAVs (0\% penetration) can easily outperform the TLC case; in Fig. \ref{678}(d) where $\lambda = 800$ veh/(hour$\cdot$lane), all the vehicles have to be CAVs (100\% penetration) in order to outperform the TLC case. Generally, higher CAV penetration is required to match the performance under TLC as traffic flow rate increases. The energy impact of traffic flow rates will be discussed in more details in Sec. \ref{sim:b}.

An additional important observation is that energy performance is not always monotonically increasing with the penetration rate value. In Fig. \ref{700}, the energy consumption under 100\% CAV
penetration is actually worse than that with 90\% and even
80\% penetration rates. This is attributed to the overly conservative nature
of our approach for determining the terminal time sequence, specifically the
fact that only one vehicle is allowed inside the MZ at any time for vehicles
traveling from different directions \eqref{eq:lateral}. On the other hand, the
collision avoidance approach adopted under mixed traffic conditions (i.e.,
\textit{Conflict Areas}) makes better use of
the MZ by allowing vehicles to share it at the same time. Such more
efficient MZ utilization reduces unnecessary travel delays.

\subsection{Energy Impact of CAV Penetration Under Different Traffic Flow Rates}

\label{sim:b}

For a dynamic system such as a transportation network, the
performance should be measured while system stability is ensured, i.e., when the intersection is
not saturated. When the intersection is saturated, it
holds too many vehicles beyond its capacity and congestion would occur. In that case, it is not possible to apply any control except to use traffic signaling.
%

The \textit{saturation flow rate} is an important concept
for evaluating the performance of transportation systems,  defined as the
headway in seconds between vehicles moving at steady state. From the
perspective of queueing theory, the intersection is a M/G/1 queue, where the
MZ is the server and the vehicles are the clients. In order for the
intersection to be stable, the vehicle arrival rate
should be less than the MZ service rate, i.e., 
$
\sum\lambda_{i}  < \mu,
$ 
where $\lambda_i$ is the arrival rate and $i$ is the index of different road segments, $\mu$ is the service rate of the MZ.
According to the recursive structure of the terminal times in \eqref{def:tf}, for vehicles
traveling on opposite roads, since they will not generate any collision
inside the MZ, they are allowed to cross the MZ at the same time. Hence, we
only need $\sum\lambda_{i}  < 2\cdot\mu$
to hold for the intersection to maintain stability. Assuming four symmetric road segments, i.e., $\lambda_i = \lambda, i \in \{1, 2, 3, 4\}$, the
saturation flow rate can be roughly estimated as $\lambda_s = 900$ veh/(hour$\cdot$lane).

The \textit{volume-to-capacity ratio}, also known as the degree of saturation, can be
calculated by dividing the actual traffic flow rate by the saturation flow rate.
Generally, an intersection with a degree of saturation less than 0.85 is considered
under-saturated and typically has sufficient capacity to maintain a stable operation.
The closer the degree of saturation is to 1, the more
sensitive the system to the natural perturbation in traffic flow. When the
degree of saturation exceeds 1, the intersection is described as over-saturated
and there is nothing we can do except to use traffic signaling. Hence, to achieve the best performance,
the traffic flow rate is set approximately under $\lambda_c $= 750 veh/(hour$\cdot$lane), where $\lambda_c$ is referred as the \textit{critical flow rate}.

\begin{figure}[tbh]
\centering
\includegraphics[width=1\columnwidth]{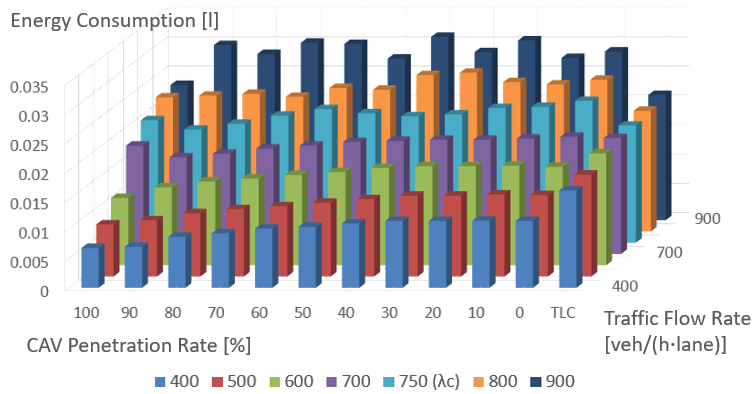}
\caption{Average energy consumption with respect to both different
traffic flow rates and CAV penetration rates.}%
\label{complete}%
\vspace{-2mm}
\end{figure}

To explore the energy impact of CAV penetration with different traffic flow rates, a comparison is presented in Fig. \ref{complete} that shows the energy consumption with respect to both different CAV penetration rates and traffic flow rates. Observe that with lower traffic flow rates, that is, when the intersection is under-saturated (i.e.,  $\lambda < \lambda_c$), the benefit obtained from CAV penetration is more significant. Even the case with no CAVs can outperform the TLC case. 
With higher traffic flow rates, that is, when the degree of saturation of the intersection is near or over 1 (i.e., $\lambda > \lambda_c$), energy consumption can hardly gain any benefit from CAV penetration. In that case, even 100\% CAV penetration cannot match the performance under TLC. This is
consistent with our expectation: when the traffic is light, the red lights
prevent some vehicles from crossing the intersection even if there is no other traffic that could generate collision inside the MZ; when the traffic is heavy, both CAVs and non-CAVs need to slow down or even stop to yield when they approach the MZ without traffic signaling and accelerate after they leave the intersection,
which may consume more energy compared to the TLC case, where
some vehicles do not need to stop during the green light phase. 

%
%

\begin{figure}[tbh]
\centering
\includegraphics[width=1\columnwidth]{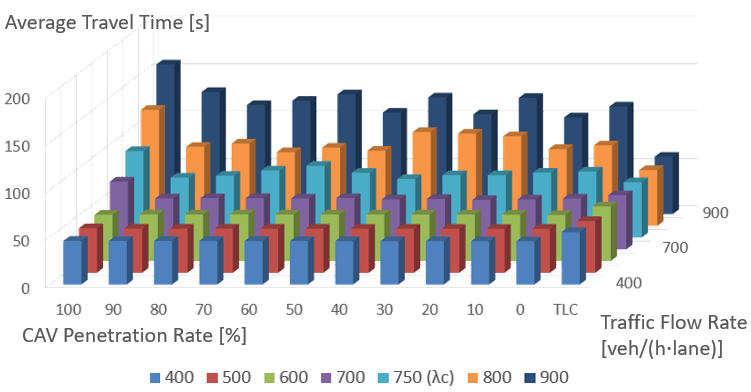}
\caption{Average travel time with respect to
different traffic flow rates and CAV penetration rates.}%
\label{complete2}%
\vspace{-2mm}
\end{figure}

The corresponding average travel times are shown in
Fig. \ref{complete2}. Observe that before the
traffic flow rate reaches the critical flow rate (i.e., $\lambda < \lambda_c$), the TLC cases are slightly outperformed in terms of travel times. This indicates that the non-signalized coordination policy (i.e., \textit{Conflict Areas}) is more effective in terms of reducing travel delay compared to TLC. This may be due to the fact that the red lights prevent some vehicles from crossing the MZ even if there is no other traffic that could generate collision, while the non-signalized coordination policy makes better use of the MZ by allowing more vehicles sharing the MZ at the same time and hence reduces travel delay. When the traffic is heavy (i.e., $\lambda > \lambda_c$), almost all the vehicles have to slow down or even stop to yield when approaching the MZ without traffic signaling,
which greatly increases the travel time.

Overall, when the traffic is light, both the energy economy and the travel time can benefit from CAV penetration. 

%

\subsection{Energy Impact of CAV Penetration Under Different Modeling Approaches}

\label{sim:c}

\begin{figure}[tbh]
\centering
\includegraphics[width=1\columnwidth]{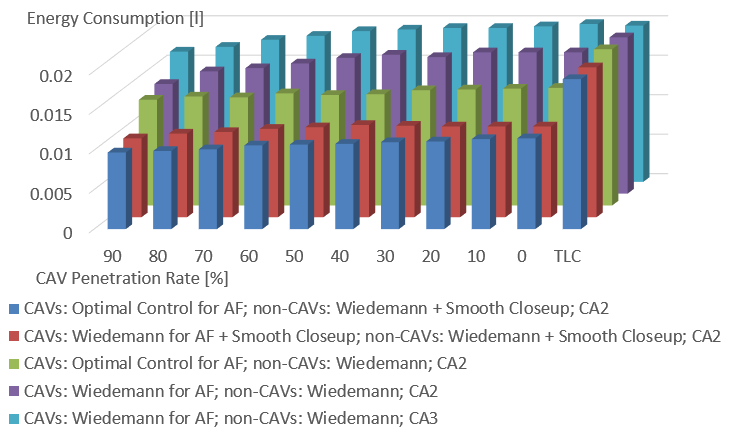} \caption{Energy
consumption under different modeling approaches.}%
\label{model_comparison}%
\vspace{-2mm}
\end{figure}

In addition to what has been discussed above, energy consumption
can also be affected by how we model the vehicle behavior and the collision
avoidance approach inside the MZ. Fig. \ref{model_comparison} shows the energy consumption under different modeling approaches:

\begin{itemize}
\item[1] For CAVs, the Wiedemann model is adopted for the AF
mode (light blue bars), while for the collision avoidance method inside the MZ,
we assume a fully controlled approach (CA3 in Fig. \ref{conflictarea}).
Observe that the energy efficiency is improving as the CAV penetration rate
increases. The benefit obtained from CAV penetration is consistent with what
we have discussed in Sec. \ref{sim:a}.

\item[2] Based on scenario 1, we modify the collision avoidance method
inside the MZ to a partially controlled approach (CA2 in Fig.
\ref{conflictarea}), where some of the vehicles do not
need to yield and hence, avoid unnecessary stops. Note that the stop-and-go process is
one of the major reasons leading to extra energy consumption. The energy
efficiency is improved  (purple bars) compared with that under scenario 1, which
indicates that more intelligent collision avoidance may save energy from
reducing unnecessary stops and travel delay.

\item[3] Based on scenario 2, we apply the optimal control in terms of
minimizing the energy consumption when CAVs are in the AF mode,
while maintaining a minimum safety following distance $\delta$ with the
preceding non-CAV if it exists. Observe that the energy performance under this case (green bars) outperforms
that under scenario 2. This validates the
effectiveness of our optimal control framework for CAVs. 

\item[4] Based on scenario 2, we adopt a less aggressive Wiedemann
car-following model by enabling the \textit{Smooth-Closeup} option. This
option forces the vehicle to decelerate in a gentle way when approaching the
preceding vehicle, so as to reduce energy consumption and passenger
discomfort, i.e., the jerk. As shown in the simulation results (red bars), the energy
consumption decreases when \textit{Smooth-Closeup} option is enabled.

\item[5] Based on scenario 4, we apply the optimal control in terms of
minimizing the energy consumption when CAVs are in the AF mode,
while maintaining a minimum safety following distance $\delta$ with the
preceding non-CAV if it exists. Compared with scenario 4, the energy
efficiency improves (dark blue bars) but the margin is not very significant. We may
reach the conclusion that the Wiedemann car-following model with the
\textit{Smooth-Closeup} option enabled has a similar impact on energy economy as the optimal control for AF mode. This is intuitively correct in
the sense that, the smoother the trip is, the less energy is consumed.
\end{itemize}


\section{Concluding Remarks}

\label{sec:con}

Earlier work \cite{Zhang2016} has established a decentralized optimal control
framework for optimally controlling CAVs crossing two adjacent intersections
in an urban area. In this paper, we extended the solution of this problem to
accommodate non-CAVs
by formulating another optimal control problem to adaptively follow the
preceding non-CAV. In addition, we investigate the energy impact of CAV
penetration under different traffic scenarios. The simulation results validate the effectiveness
of CAV penetration, and as the CAV
penetration rate increases, the benefit becomes more significant. Such results
can be consistently observed across a range of different traffic intensity
settings as long as the traffic flow rate is below the critical flow rate; this provides strong evidence of the advantages of incorporating CAVs into current traffic systems.

Ongoing research is considering turns (see \cite{Zhang2017}) and
lane changing in the intersection with a diverse set of CAVs and exploring the
associated tradeoffs between the intersection throughput and energy
consumption of each individual vehicle. Future research should also
investigate the potential to further maximize the traffic throughput by re-evaluation of the vehicle crossing sequence as new vehicles enter the CZ.

\bibliographystyle{IEEETran}
\bibliography{CCTA2018}

\begin{thebibliography}{10}
\providecommand{\url}[1]{#1}
\csname url@rmstyle\endcsname
\providecommand{\newblock}{\relax}
\providecommand{\bibinfo}[2]{#2}
\providecommand\BIBentrySTDinterwordspacing{\spaceskip=0pt\relax}
\providecommand\BIBentryALTinterwordstretchfactor{4}
\providecommand\BIBentryALTinterwordspacing{\spaceskip=\fontdimen2\font plus
\BIBentryALTinterwordstretchfactor\fontdimen3\font minus
  \fontdimen4\font\relax}
\providecommand\BIBforeignlanguage[2]{{%
\expandafter\ifx\csname l@#1\endcsname\relax
\typeout{** WARNING: IEEEtran.bst: No hyphenation pattern has been}%
\typeout{** loaded for the language `#1'. Using the pattern for}%
\typeout{** the default language instead.}%
\else
\language=\csname l@#1\endcsname
\fi
#2}}

\bibitem{Fleck2015}
J.~L. Fleck, C.~G. Cassandras, and Y.~Geng, ``Adaptive quasi-dynamic traffic
  light control,'' \emph{IEEE Transactions on Control Systems Technology},
  2015, DOI: 10.1109/TCST.2015.2468181, to appear.

\bibitem{Levine1966}
W.~Levine and M.~Athans, ``{On the optimal error regulation of a string of
  moving vehicles},'' \emph{IEEE Transactions on Automatic Control}, vol.~11,
  no.~3, pp. 355--361, 1966.

\bibitem{Dresner2004}
K.~Dresner and P.~Stone, ``{Multiagent traffic management: a reservation-based
  intersection control mechanism},'' in \emph{Proceedings of the Third
  International Joint Conference on Autonomous Agents and Multiagents Systems},
  2004, pp. 530--537.

\bibitem{Dresner2008}
------, ``{A Multiagent Approach to Autonomous Intersection Management},''
  \emph{Journal of Artificial Intelligence Research}, vol.~31, pp. 591--653,
  2008.

\bibitem{DeLaFortelle2010}
A.~{de La Fortelle}, ``{Analysis of reservation algorithms for cooperative
  planning at intersections},'' \emph{13th International IEEE Conference on
  Intelligent Transportation Systems}, pp. 445--449, Sept. 2010.

\bibitem{Huang2012}
S.~Huang, A.~Sadek, and Y.~Zhao, ``{Assessing the Mobility and Environmental
  Benefits of Reservation-Based Intelligent Intersections Using an Integrated
  Simulator},'' \emph{IEEE Transactions on Intelligent Transportation Systems},
  vol.~13, no.~3, pp. 1201,1214, 2012.

\bibitem{Li2006}
L.~Li and F.-Y. Wang, ``{Cooperative Driving at Blind Crossings Using
  Intervehicle Communication},'' \emph{IEEE Transactions in Vehicular
  Technology}, vol.~55, no.~6, pp. 1712,1724, 2006.

\bibitem{Yan2009}
F.~Yan, M.~Dridi, and A.~{El Moudni}, ``{Autonomous vehicle sequencing
  algorithm at isolated intersections},'' \emph{2009 12th International IEEE
  Conference on Intelligent Transportation Systems}, pp. 1--6, 2009.

\bibitem{Zohdy2012}
I.~H. Zohdy, R.~K. Kamalanathsharma, and H.~Rakha, ``{Intersection management
  for autonomous vehicles using iCACC},'' \emph{2012 15th International IEEE
  Conference on Intelligent Transportation Systems}, pp. 1109--1114, 2012.

\bibitem{Zhu2015}
F.~Zhu and S.~V. Ukkusuri, ``{A linear programming formulation for autonomous
  intersection control within a dynamic traffic assignment and connected
  vehicle environment},'' \emph{Transportation Research Part C: Emerging
  Technologies}, 2015.

\bibitem{Lee2012}
J.~Lee and B.~Park, ``{Development and Evaluation of a Cooperative Vehicle
  Intersection Control Algorithm Under the Connected Vehicles Environment},''
  \emph{IEEE Transactions on Intelligent Transportation Systems}, vol.~13,
  no.~1, pp. 81--90, 2012.

\bibitem{Miculescu2014}
D.~Miculescu and S.~Karaman, ``{Polling-Systems-Based Control of
  High-Performance Provably-Safe Autonomous Intersections},'' in \emph{53rd
  IEEE Conference on Decision and Control}, 2014.

\bibitem{Rios-Torres}
J.~Rios-Torres and A.~A. Malikopoulos, ``{A Survey on the Coordination of
  Connected and Automated Vehicles at Intersections and Merging at Highway
  On-Ramps},'' \emph{IEEE Transactions on Intelligent Transportation Systems},
  2016 (forthcoming).

\bibitem{ZhangMalikopoulosCassandras2016}
Y.~Zhang, A.~A. Malikopoulos, and C.~G. Cassandras, ``Optimal control and
  coordination of connected and automated vehicles at urban traffic
  intersections,'' in \emph{Proceedings of the American Control Conference},
  2016, pp. 6227--6232.

\bibitem{gilbert1976vehicle}
E.~G. Gilbert, ``Vehicle cruise: Improved fuel economy by periodic control,''
  \emph{Automatica}, vol.~12, no.~2, pp. 159--166, 1976.

\bibitem{hooker1988optimal}
J.~Hooker, ``Optimal driving for single-vehicle fuel economy,''
  \emph{Transportation Research Part A: General}, vol.~22, no.~3, pp. 183--201,
  1988.

\bibitem{hellstrom2010design}
E.~Hellstr{\"o}m, J.~{\AA}slund, and L.~Nielsen, ``Design of an efficient
  algorithm for fuel-optimal look-ahead control,'' \emph{Control Engineering
  Practice}, vol.~18, no.~11, pp. 1318--1327, 2010.

\bibitem{li2012minimum}
S.~E. Li, H.~Peng, K.~Li, and J.~Wang, ``Minimum fuel control strategy in
  automated car-following scenarios,'' \emph{IEEE Transactions on Vehicular
  Technology}, vol.~61, no.~3, pp. 998--1007, 2012.

\bibitem{dresner2007sharing}
K.~M. Dresner and P.~Stone, ``Sharing the road: Autonomous vehicles meet human
  drivers.'' in \emph{IJCAI}, vol.~7, 2007, pp. 1263--1268.

\bibitem{guler2014using}
S.~I. Guler, M.~Menendez, and L.~Meier, ``Using connected vehicle technology to
  improve the efficiency of intersections,'' \emph{Transportation Research Part
  C: Emerging Technologies}, vol.~46, pp. 121--131, 2014.

\bibitem{Malikopoulos2016}
A.~A. Malikopoulos, C.~G. Cassandras, and Y.~Zhang, ``A decentralized
  energy-optimal control framework for connected automated vehicles at
  signal-free intersections,'' \emph{Automatica}, 2017 (to appear).

\bibitem{Zhang2016}
Y.~Zhang, C.~G. Cassandras, and A.~A. Malikopoulos, ``Optimal control of
  connected automated vehicles at urban traffic intersections: A feasibility
  enforcement analysis,'' in \emph{Proceedings of the 2017 American Control
  Conference}, 2017, pp. 3548--3553.

\bibitem{wiedemann1974simulation}
R.~Wiedemann, ``Simulation des strassenverkehrsflusses.'' 1974.

\bibitem{Zhang2017}
Y.~Zhang, A.~A. Malikopoulos, and C.~G. Cassandras, ``Decentralized optimal
  control for connected automated vehicles at intersections including left and
  right turns,'' in \emph{56th IEEE Conference on Decision and Control}, 2017,
  pp. 4428--4433.

\end{thebibliography}

\end{document}